\documentclass{article}
\newcommand{\proof}{\noindent {\bf Proof: }}

\newtheorem{theorem}{Theorem}

\newtheorem{corollary}{Corollary}
\newtheorem{lemma}{Lemma}

\def\qed{\hfill $\Box$}

\usepackage[]{amssymb}
\usepackage{graphicx}

\begin{document}
\title{Normally distributed probability measure on the metric space of norms}
\author{\'A.G.Horv\'ath}
\date{2011 May}

\maketitle

\begin{abstract}
In this paper we propose a method to construct probability measures on the space of convex bodies with a given pushforward distribution. Concretely we show that there is a measure on the metric space of centrally symmetric convex bodies, which pushforward by the thinness mapping produces a probability measure of truncated normal distribution on the interval of its range. Improving the construction we give another (more complicated) one with the following additional properties; the neighborhoods have positive measure, the set of polytopes has zero measure and the set of smooth bodies has measure 1, respectively.
\end{abstract}



\section{Introduction}

In this paper we shall investigate the probability space of norms defined on a real, $n$-dimensional Euclidean space $V$. A norm function on $V$ defined by its unit ball $K$, which is (in a fixed, cartesian coordinate system of the Euclidean vector space $(V, \langle \cdot, \cdot \rangle)$ with origin $O$ is a centrally symmetric in $O$) convex body. Such bodies give a closed proper subset $\mathcal{K}_0$ \footnote{We rather denote in this paper the space of $O$-symmetric convex bodies by $\mathcal{K}_0$ as the space of convex bodies with centroid $O$.} of the space of convex bodies $\mathcal{K}$ of $(V, \langle \cdot, \cdot \rangle)$. It is known that the Hausdorff distance (denoted by $\delta^h$) is a metric on $\mathcal{K}$ and with this metric $(\mathcal{K},\delta^{h})$ is a locally compact space. (See in \cite{gruber 1},\cite{gruber 2}.) Thus there should be many measures available on these space. Unfortunatelly this is not so. Bandt and Baraki in \cite{bandt} proved answering to a problem of McMullen \cite{schneider} that there is no positive $\sigma$-finite Borel measure on it which is invariant with respect to all isometries of $(\mathcal{K},\delta^h)$ into itself. This result exclude the possibility of the existence of a natural volume-type measure. It was a natural question that can whether be found such a $\sigma $-finite Borel measure on $\mathcal{K}$ which holds the property that it is non-zero for any open set of $\mathcal{K}$ and invariant under rigid motions of the embedding vector space. This long standing question was answered in the last close by Hoffmann in \cite{hoffmann}. His result can be summarized as follows.
Each $\sigma$-finite rotation and translation
invariant Borel measure on $(\mathcal{K},\delta^h)$ is the vague limit of such measures and that each $\sigma $-finite
Borel measure on $(\mathcal{K},\delta^h)$ is the vague limit of measures of the form
$$
\sum\limits_{i=1}^{\infty} \alpha_n\delta_{K_n},
$$
where $\{K_n \mbox{ , } n\in \mathbb{N}\}$ is a countable, dense subset of $(\mathcal{K},\delta^h)$, $(\alpha_n)$ is a sequence of positive real numbers for which $\sum \limits_{i=1}^\infty \alpha_n <\infty $ and $\delta_{K_n}$ denote the Dirac measure concentrated at $K_n$.

Hoffmann also observed that a result of B\'ar\'any \cite{barany} "suggest that it might not be possible to define a "uniform" probability measure on the set of all polytopes which have rational vertices and are contained in the unit ball". The known concept of Gaussian random convex bodies \cite{molchanov} gives a poor class of Gaussian measures because of a random convex body is Gaussian if and only if there exists a deterministic body and a Gaussian random vector such that the random body is the sum of the deterministic one and the random vector almost surely. He asked "whether there exists an alternative approach to "Gaussian" random convex bodies which yields a richer class of "Gaussian" measures on $(\mathcal{K},\delta^h)$.

In present paper we propose that define "Gaussian" (or other type) probability measure on $(\mathcal{K},\delta^h)$ with respect to a given measurable function of the space. Our observation is that on certain probability space the uniformity or normality properties could be only "relative" one and thus we can require these properties in their impacts through a given function of the space. More precisely, we require the normality or uniformity on a pushforward measure by a given geometric function of the elements of the space (here on the space of convex bodies). To this purpose we will use the thinness function $\alpha_0(K)$ of $K$ defined by the help of the concepts of diameter $d(K)$ and width $w(K)$, respectively, by the equality:
$$
\alpha_0(K)=\frac{d(K)}{w(K)+d(K)}.
$$
As a concrete construction we will give a probability measure on  $(\mathcal{K}_0,\delta^h)$ which pushforward measure by the function $\alpha_0(K)$ has truncated normal distribution on the range interval $[\frac{1}{2},1)$ (Theorem 2). We note that a method which sends a convex body to a well-defined $O$-symmetric convex body by a continuous mapping, define a pullback measure of to the all space of convex bodies. The pushforward of this pullback measure by the composition of the mapping of the method and the function $\alpha_0(K)$ has the same properties as the measure of $(\mathcal{K}_0,\delta^h)$. To this purpose we can use the Minkowski symmetrization process sending a body $K$ into the body $\frac{1}{2}(K+(-K))$ with the same diameter, width and thinness (Corollary 1). Our last statement (Theorem 3) that the previously construction can be modified such that the set of basic bodies will be dense and countable with smooth elements. Thus the set of polytopes has zero measure, the set of smooth bodies has measure 1, and every neighborhood has positive measure.

\section{The thinness function}

Let denote by $\alpha_0(K)$ the number
$$
\alpha_0(K)=\frac{d(K)}{w(K)+d(K)},
$$
where $w(K)$ and $d(K)$ denote the width and diameter of the body $K$, respectively. This number characterize the thinness of the convex body $K$, it is $\frac{1}{2}$ in the case of the Euclidean ball only and is equal to $1$ if $K$ has of dimension less or equal to $n-1$.

Let now $B_E$ be the unit ball of the embedding Euclidean space and let define the unit sphere of $\mathcal{K}_{0}$ around $B_E$ by the equality:
$$
\mathcal{K}_{0}^1:=\{ K\in \mathcal{K}_{0} \mid \delta^h(K,B_E)=1\}.
$$

The following lemma shows the usable of the thinness function in our investigation.

\begin{lemma}
If $K\in \mathcal{K}_{0}^1$ and $\alpha_0:=\alpha_0(K)$ is the thinness of $K$ then we have
$$
\delta^h(\alpha K,B_E)=\left\{ \begin{array}{lcc}
2\alpha-1 & \mbox{ if } & \alpha_0 \leq \alpha\\
2\alpha+1-2\frac{\alpha}{\alpha_0} & \mbox{ if } & 0\leq \alpha <\alpha_0.
\end{array} \right.
$$
\end{lemma}

\proof
Assume that $\delta^h(K,B_E)$ is the distance of the points $x\in \mathrm{bd} B_E$ and $y\in \mathrm{ bd} K$. Then $\|y\|_E=\|x\|_E+1=2$ and $0,x,y$ are collinear. (We note that the norm of the point $y$ is also the half of the diameter $d(K)$ of $K$ with respect to the Euclidean metric.) This implies that for $\alpha >1$ the points $\frac{1}{\alpha}x$ and $y$ give a segment with length $\delta^h\left(K,\frac{1}{\alpha}B_E\right)$ and thus
$$
\delta^h\left(K,\frac{1}{\alpha}B_E\right)=\|y-\frac{1}{\alpha}x\|_E=\|y\|_E-\frac{1}{\alpha}\|x\|_E=2-\frac{1}{\alpha}
$$
holds. If $\alpha<1$ then the situation is a little bit more complicated. In this case there is a real number $\alpha _0\in \left[\frac{1}{2},1\right) $ such that if $\alpha_0\leq \alpha <1$ then again
$$
\delta^h\left(K,\frac{1}{\alpha}B_E\right)=\|y-\frac{1}{\alpha}x\|_E=\|y\|_E-\frac{1}{\alpha}\|x\|_E=2-\frac{1}{\alpha}
$$
but for $\alpha_0\geq \alpha >0$ we have a new pair of points  $y'\in \mathrm {bd} K$ and $x'\in \mathrm {bd} B_E$  where the distance attained. The point $y'$ is a point of $ \mathrm {bd} K$ with minimal norm and we have the equality
$$
\frac{1}{\alpha_0}-\|y'\|=2-\frac{1}{\alpha_0}.
$$
Thus the norm of $y'$ is equal to $2(\frac{1}{\alpha_0}-1)$. In this case
$$
\delta^h\left(K,\frac{1}{\alpha}B_E\right)=\left\|-y'+\frac{1}{\alpha}x'\right\|_E=\frac{1}{\alpha}-2\left(\frac{1}{\alpha_0}-1\right)=2+\frac{1}{\alpha}-\frac{2}{\alpha_0}.
$$
We thus have the equality
$$
\delta^h(\alpha K,B_E)=\alpha \delta^h\left(K,\frac{1}{\alpha}B_E\right)=\left\{ \begin{array}{lcc}
2\alpha-1 & \mbox{ if } & \alpha_0 \leq \alpha\\
2\alpha+1-2\frac{\alpha}{\alpha_0} & \mbox{ if } & 0\leq \alpha <\alpha_0.
\end{array} \right.
$$
The constant $\alpha_0$ depends only on the body $K$ and it has the following geometric meaning. $\|y'\|_E=\frac{2}{\alpha_0}-2$ is the half of the width $w(K)$ of the centrally symmetric body $K$, because it is a point on $\mathrm{bd }K$ with minimal norm. So we can see that
$$
\frac{1}{2}\leq \alpha_0=\frac{2}{\|y'\|_E+2}=\frac{d(K)}{w(K)+d(K)}<1
$$
as we stated.
\qed

\section{Measure on $\mathcal{K}_{0}^1$ with uniform pushforward.}

 We now construct a measure on $\mathcal{K}_0^1$ which pushforward by the thinness function has uniform distribution. To this (following Hoffmann's paper) we introduce the orbits of a body $K$ about the special orthogonal group $SO(n)$ by $[K]$. These are compact subsets of $\mathcal{K}_0^1$, and if we consider an open subset of $\mathcal{K}_0^1$ then the union of the corresponding orbits is also open. Hence there exists a measurable mapping $s:\mathcal{K}_0^1\longrightarrow \mathcal{K}_0^1$ such that $s(K)=s(K')$ if and only if $K$ and $K'$ are on the same orbit. Let $\widetilde{\mathcal{K}_0^1}:=\{ K\in \mathcal{K}_0^1 \mbox{ , } s(K)=K\}$ which is measurable subset of $\mathcal{K}_0^1$. We equip it with the induced topology of $\mathcal{K}_0^1$. Finally let $\Phi_{2a}^1:\widetilde{\mathcal{K}_0^1}\times SO(n)\longrightarrow \mathcal{K}_0^1$ is the mapping defined by the equality:
$$
\Phi_{2a}^1(K,\Theta)=\Theta K.
$$
Our notation is analogous with the notation of \cite{hoffman}. It was proved in \cite{hoffmann} (Lemma 2) that a non-trivial $\sigma$-finite measure $\mu_0$ on $\mathcal{K}_0$ is invariant under rotations (meaning that for $\Theta \in SO(n)$ we have $\mu_0(\mathcal{A})=\mu(\Theta \mathcal{A})$ for all Borel sets $\mathcal{A}$ of $\mathcal{K}_0$) if and only if there exists a $\sigma$-finite measure $\widetilde{\mu_0}$ on $\widetilde{\mathcal{K}_0}$ such that $\mu_0=\Phi_{2a}(\widetilde{\mu_0}\otimes \nu_n)$, where $\nu_n$ is the Haar measure on $SO(n)$. It is obvious that in the case of $\mathcal{K}_0^1$ there is a similar result by our mapping $\Phi_{2a}^1(K,\Theta)$ which is the restriction of Hoffmann's map $\Phi_{2a}(K,\Theta)$ onto the set  $\mathcal{K}_0^1$.

First choose a countable system of bodies $K_m$ to define a probability measure on $\widetilde{\mathcal{K}_0^1}$. Without loss of generality we may assume that each of the bodies of $\widetilde{\mathcal{K}_0^1}$ has a common diameter of length $4$ denoted by $d$, which lies on the $n^{th}$ axe of coordinates (hence it is the convex hull of the points $\{2e_n,-2e_n\}$). Consider the set of diadic rational numbers in $(0,2]$. We can write them as follows:
$$
\left\{ m(n,k):=\frac{k}{2^n} \mbox{ where } n=0, \cdots \infty \mbox{ and for a fixed } n \mbox{, }0<k\leq 2^{n+1}\right\}.
$$
 Define the body $K_{m(n,k)}$ as the convex hull of the union of the segment $d$ and the ball around the origin with radius $m(n,k)$. For each $n$ we have $2^{n+1}$ such bodies, thus the definition
$$
\widetilde{\mu_0^1}:=\lim\limits_{n\rightarrow \infty}\sum\limits_{k=1}^{2^{n+1}}\frac{1}{2^{n+1}}\delta_{K_{m(n,k)}}
$$
define a probability measure on $\widetilde{\mathcal{K}_0^1}$. (The limit is the vague limit (or limit with respect to weak convergence) of measures.) In fact,
$$
\widetilde{\mu_0^1}\left(\widetilde{\mathcal{K}_0^1}\right)=\lim\limits_{n\rightarrow \infty}\sum\limits_{k=0}^{2^{n+1}}\frac{1}{2^{n+1}}\delta_{K_{m(n,k)}}\left(\widetilde{\mathcal{K}_0^1}\right)=1
$$
\begin{lemma}
The pushforward measure
$w(K)^{-1}(\widetilde{\mu_0^1})$ has uniform distribution on the interval $(0,4]$. ( $w(K)$ means the width of the body $K$.)
\end{lemma}

\proof Let $B'=(0,x]$ be a level set of  $(0,4]$. By definition
$$
w(K)^{-1}(\widetilde{\mu_0^1})(B')=\widetilde{\mu_0^1}\left(\left\{K\in \widetilde{\mathcal{K}_0^1} \mid w(K)\in B'\right\}\right)=
$$
$$
=\lim\limits_{n\rightarrow \infty}\sum\limits_{K_{m(n,k)}\in w(K)^{-1}(B') \atop 0< k\leq 2^{n+1}}\frac{1}{2^{n+1}}=
\lim\limits_{n\rightarrow \infty}\sum\limits_{2m(n,k)\in B'}\frac{1}{2^{n+1}}=
$$
$$
=\lim\limits_{n\rightarrow \infty}\sum\limits_{2m(n,k)<x}\frac{1}{2^{n+1}}=\lim\limits_{n\rightarrow \infty}\sum\limits_{k=1}^{2^{n-1}x}\frac{1}{2^{n+1}}= \frac{x}{4}
$$
showing that $w(K)^{-1}(\widetilde{\mu_0^1})$ is the uniform distribution of the interval $(0,4]$.
\qed

The Gaussian measure $\gamma $ of the $n^2$-dimensional matrix space $\mathbb{R}^{n\times n}$ defined by the density function $G(X)$
$$
G(X)\mathrm{d\lambda^{n^2}}:=\frac{1}{\left(\sqrt{2\pi}\right)^{n^2}}e^{-\frac{1}{2}\mathrm{Tr}(X^TX)}\mathrm{ d\lambda^{n^2} },
$$
where $\mathrm{ d\lambda^{n^2} }$ is the $n^2$-dimensional Lebesgue measure. The Haar measure $\nu_n$ is the pushforward measure of the Gaussian measure by the mapping $M$ defined by the Gram-Schmidt process (see in \cite{lee}). In fact, if $GL(n,\mathbb{R})$ is the group of nonsingular matrices then $M$ is a mapping from
$GL(n,\mathbb{R})$ to $O(n)$ with the following properties:
\begin{enumerate}
\item surjective;
\item if $B$ is a Borel set then we have $M(QB)=QM(B)$ for every $Q\in O(n)$;
\item  $M^{-1}(QB)=QM^{-1}(B)$.
\end{enumerate}
It can be proved that for a Borel set $B$ of $O(n)$ we have
$$
\nu_n(B)=\gamma\left(M^{-1}(B)\right).
$$
Furthermore $\nu_n$ is a probability measure because it can be seen that a matrix invertible (with respect to the Gaussian measure) with probability 1 and thus
$$
\nu_n(O(n))=\gamma(GL(n,\mathbb{R}))=1.
$$
Since Haar measure by definition invariant under orthogonal transformations it is the unique "uniform" (geometric volume) distribution on $O(n)$ and thus on $SO(n)$, too.

We now state the following:
\begin{theorem}
Let define the measure $\widetilde{\nu_0^1}$ by density function $\mathrm{d \widetilde{\nu_0^1}}=\frac{4}{(w+4)^2}\mathrm{d \widetilde{\mu_0^1}}$. Then
$$
\alpha_0(K)^{-1}\left(\Phi_{2a}^1\left(\widetilde{\nu_0^1}\otimes \nu_n\right)\right)
$$
is a probability measure with uniform distribution on $[\frac{1}{2},1)$.
\end{theorem}

\proof We are stating that the pushforward measure
$$
\alpha_0(K)^{-1}\left(\Phi_{2a}^1\left(\left(\widetilde{\nu_0^1}\otimes \nu_n\right)\right)\right)
$$
has uniform distribution on $[\frac{1}{2}, 1)$
if and only if the pushforward measure
$$
w(K)^{-1}\left(\widetilde{\mu_0^1}\right)
$$
has uniform distribution on $(0,4]$.
To prove this consider a Borel set $B$ of $[\frac{1}{2},1)$ and its image $B'$ under the bijective transformation
$$
\tau:t\mapsto \tau(t):=\frac{4}{t}-4.
$$
Of course $B'$ is a Borel set of the interval $(0,4]$ which is the image of $[\frac{1}{2},1)$ with respect to $\tau$.
We now have that
$$
\int\limits_B\mathrm{ d }\alpha_0(K)^{-1}\left(\Phi_{2a}^1\left(\widetilde{\nu_0^1}\otimes \nu_n\right)\right)=\alpha_0(K)^{-1}\left(\Phi_{2a}^1\left(\widetilde{\nu_0^1}\otimes \nu_n\right)\right)(B)=
$$
$$
=\Phi_{2a}^1\left(\widetilde{\nu_0^1}\otimes \nu_n\right)(\alpha_0(K)^{-1}(B))=
$$
$$
=\widetilde{\nu_0^1}\left(\left(\Phi_{2a}^1\right)^{-1}_1\left((\alpha_0(K)^{-1}(B))\right)\right)\nu_n\left(\left(\Phi_{2a}^1\right)^{-1}_2\left( \alpha_0(K)^{-1}(B))\right)\right)
$$
where $ \left(\Phi_{2a}^1\right)^{-1}_1$ and $\left(\Phi_{2a}^1\right)^{-1}_2$ means the components of the set-valued inverse of the function $\Phi_{2a}^1$, respectively. Since $\left(\Phi_{2a}^1\right)^{-1}_2\left( \alpha_0(K)^{-1}(B))\right)$ is the group $O(n)$ we have that
$$
\int\limits_B\mathrm{ d }\alpha_0(K)^{-1}\left(\Phi_{2a}^1\left(\widetilde{\nu_0^1}\otimes \nu_n\right)\right)=\widetilde{\nu_0^1}\left(\left(\Phi_{2a}^1\right)^{-1}_1\left(\alpha_0(K)^{-1}(B)\right)\right)=
$$
$$
=\int\limits_{\left(\Phi_{2a}^1\right)^{-1}_1\left(\alpha_0^{-1}(B)\right)}\mathrm{ d }\widetilde{\nu_0^1}.
$$
On the other hand
$$
\left(\Phi_{2a}^1\right)^{-1}_1\left(\alpha_0^{-1}(B)\right)=\left\{\widetilde{K}\in \widetilde{\mathcal{K}_0^1} \mid \alpha_0(\widetilde{K})=\frac{4}{w(\widetilde{K})+4}\in B\right\}=
$$
$$
=\left\{\widetilde{K}\in \widetilde{\mathcal{K}_0^1}  \mid w(\widetilde{K})\in B'=\frac{4}{B}-4\right\}
$$
implying that
$$
\int\limits_{\left(\Phi_{2a}^1\right)^{-1}_1\left(\alpha_0^{-1}(B)\right)}\mathrm{ d }\widetilde{\nu_0^1}=\int\limits_{\left\{\widetilde{K}\in \widetilde{\mathcal{K}_0^1}  \mid w(\widetilde{K})\in B'\right\}}\frac{4}{(w+4)^2}\mathrm{ d }\widetilde{\mu_0^1},
$$
and it is equal to
$$
\int\limits_{\tau\in B'}\frac{4}{(4+\tau)^2}\mathrm{ d\tau }=\int\limits_{t\in B}\mathrm{ dt }
$$
if and only if $w(K)^{-1}\left(\widetilde{\mu_0^1}\right)$ has uniform distribution on $(0,4]$ as we stated.

Since Lemma 2 says that $w(K)^{-1}\left(\widetilde{\mu_0^1}\right)$ has uniform distribution on the interval $[0,4]$ we also proved the theorem.
\qed

Let denote by $\nu_0^1$ the measure
$$
\Phi_{2a}^1\left(\widetilde{\nu_0^1}\otimes \nu_n\right).
$$

\section{Measure on $\mathcal{K}_{0}$ with normal pushforward.}

Finally we can identify $\mathcal{K}_{0}$ with $\mathcal{K}_{0}^1\times [0,\infty)$.
To this end let $\Phi_4$ be the mapping
$$
\Phi_4:(K,\alpha) \mapsto \alpha K.
$$
\begin{lemma}
From the image $K'=\Phi_4(K)$ we can determine uniquely the body $K$ and the constant $\alpha $.
\end{lemma}

\proof
$
K'=\alpha K
$
implies that $\alpha_0(K)=\alpha_0(K')=\frac{d(K')}{w(K')+d(K')}$ and thus $\alpha_0(K)$ is uniquely determined. We also know the value of
$$
\alpha':=\delta^h(\alpha K,B_E).
$$
We consider two cases. In the first case we assume that $\alpha \geq \alpha_0$ and hence by Lemma 1 we get that
$$
\alpha'=2\alpha-1 \mbox{ or } \alpha =\frac{\alpha'+1}{2}
$$
and in the second one we assume $0\leq\alpha\leq \alpha_0$ then we have
$$
\alpha'= 2\alpha+1-2\frac{\alpha}{\alpha_0} \mbox{ or } \alpha=\frac{\alpha'-1}{2-\frac{2}{\alpha_0}}=\frac{\alpha_0(\alpha'-1)}{2(\alpha_0-1)}
$$
From these equalities we get that the first case implies
$$
\alpha_0\leq \frac{\alpha'+1}{2} \mbox{ so }  \alpha'\geq 2\alpha_0 -1
$$
and in the second one we have
$$
\alpha_0\geq \frac{\alpha_0(\alpha'-1)}{2(\alpha_0-1)}\geq 0.
$$
Hence we have
$$
2\alpha_0-1\geq \alpha' \geq 0.
$$
So first we determine $\alpha'$ and the value
$$
2\alpha_0-1=\frac{2d(K)}{w(K)+d(K)}-1=\frac{d(K)-w(K)}{d(K)+w(K)}.
$$
Then using the above equalities we can calculate $\alpha$ which is uniquely determined. Now $K$ is equal to $\frac{1}{\alpha}K'$.
\qed

Denote by $\Phi_4^{-1}(K'):=\left(\left(\Phi_4^{-1}\right)_1(K'),\left(\Phi_4^{-1}\right)_2(K')\right)$ the above determined pair $(K,\alpha)$.
If we have a $\sigma$-finite measure $\nu_0^1$ on $\mathcal{K}_{0}^1$ then we also have a $\sigma$-finite measure $\nu_0$ on $\mathcal{K}_{0}$ by the definition
$$
\nu_0=\Phi_4(\nu_0^1\otimes \nu),
$$
where $\nu$ is a $\sigma$-finite measure on $(0,\infty)$.

Let define now the set function $p(\mathcal{A})$ as follows. If the set $\mathcal{A}\subset \mathcal{K}_{0}$ $\nu_0$ measurable let be
$$
p(\mathcal{A}):=\frac{1}{\sqrt{2\pi\sigma^2}}\int\limits_{K'\in \mathcal{A}}e^{-\frac{\left(\delta^h \left(B_E,\frac{\alpha_0(K')}{\Phi_4^{-1}(K')_2}K'\right)\right)^2}{2\sigma^2}}\mathrm{d\nu_0}.
$$
The following theorem is our main result.
\begin{theorem}
If $\nu_0^1$ is such a probability measure on $\mathcal{K}_0^1$ for which $\alpha_0(K)^{-1}(\nu_0^1)$ has uniform distribution, $\nu_0=\Phi_4(\nu_0^1\otimes \nu)$ where $\nu $ is a probability measure on $(0,\infty)$ and  $\Phi$ is the probability function of the standard normal distribution then
$$
P(\mathcal{A}):=\frac{4p(\mathcal{A})}{\left(\Phi\left(\frac{1}{\sigma}\right)-\Phi(0)\right)}=
$$
$$
=\frac{4}{\left(\Phi\left(\frac{1}{\sigma}\right)-\Phi(0)\right)\sqrt{2\pi\sigma^2}}\int\limits_{K'\in \mathcal{A}}e^{-\frac{\left(\delta^h \left(B_E,\frac{\alpha_0(K')}{\Phi_4^{-1}(K')_2}K'\right)\right)^2}{2\sigma^2}}\mathrm{d\nu_0}
$$
is a probability measure on $\mathcal{K}_0$. Moreover $\alpha_0(K)^{-1}(P)$ has truncated normal distribution on the interval $[\frac{1}{2},1)$, (with mean $\frac{1}{2}$ and variance $\left(\frac{\sigma}{2}\right)^2$), so
$$
\alpha_0(K)^{-1}(P)\left(\left\{\frac{1}{2}\leq t\leq c\right\}\right)=P\left(\left\{\mathcal{K}\in \mathcal{K}_0 \mid \alpha_0(K)\leq c\right\}\right)=\frac{\Phi\left(\frac{c-\frac{1}{2}}{\frac{\sigma}{2}}\right)-\Phi(0)}{\Phi\left(\frac{1}{\sigma}\right)-\Phi(0)}.
$$
\end{theorem}

\proof
$$
p(\mathcal{A})=\frac{1}{\sqrt{2\pi\sigma^2}}\int\limits_{K\in \left(\Phi_4^{-1}\right)_1(\mathcal{A})}\int\limits_{\alpha \in \left(\Phi_4^{-1}\right)_2(\mathcal{A})}e^{-\frac{\left(\delta^h \left(B_E,\frac{\alpha_0(K')}{\alpha}\alpha K\right)\right)^2}{2\sigma^2}}\mathrm{d\nu}\mathrm{d\nu_0^1}
$$
however $\alpha_0(K')=\alpha_0(K)$ so it is equal to
$$
\frac{1}{\sqrt{2\pi\sigma^2}}\int\limits_{K\in \left(\Phi_4^{-1}\right)_1(\mathcal{A})}\left(\int\limits_{\alpha\in\left(\Phi_4^{-1}\right)_2(\mathcal{A})}e^{-\frac{\alpha_0(K)'^2}{2\sigma^2}} \mathrm{d\nu}\right)\mathrm{d\nu_0^1}=
$$
$$
=\frac{1}{\sqrt{2\pi\sigma^2}}\int\limits_{K\in \left(\Phi_4^{-1}\right)_1(\mathcal{A})}\left(\int\limits_{\alpha K\in \mathcal{A}, \atop \alpha\geq \alpha_0(K)}e^{-\frac{(2\alpha_0(K)-1)^2}{2\sigma^2}}\mathrm{d\nu}+\right.
$$
$$
\left.+\int\limits_{\alpha K\in \mathcal{A}\atop 0\leq\alpha\leq \alpha_0(K)}e^{-\frac{\left(2\alpha_0(K)+1-2\frac{\alpha_0(K)}{\alpha_0(K)}\right)^2}{2\sigma^2}}\mathrm{d\nu}\right)\mathrm{d\nu_0^1}=
$$
$$
=\frac{1}{\sqrt{2\pi\sigma^2}}\int\limits_{K\in \left(\Phi_4^{-1}\right)(\mathcal{A})_1}\left(\int\limits_{\alpha\in\left(\Phi_4^{-1}\right)_2(\mathcal{A})} e^{-\frac{(2\alpha_0(K)-1)^2}{2\sigma^2}}\mathrm{d\nu}\right)\mathrm{d\nu_0^1}=
$$
$$
=\frac{\nu\left(\alpha\in\left(\Phi_4^{-1}\right)_2(\mathcal{A})\right)} {\sqrt{2\pi\sigma^2}}\int\limits_{K\in \left(\Phi_4^{-1}\right)(\mathcal{A})_1}e^{-\frac{(2\alpha_0(K)-1)^2}{2\sigma^2}}\mathrm{d\nu_0^1}.
$$
For $\mathcal{A}=\mathcal{K}_0$ we have that it is equal to
$$
=\frac{\nu((0,\infty))}{\sqrt{2\pi\sigma^2}}\int\limits_{\frac{1}{2}}^{1}e^{-\frac{1}{2}\left(\frac{t- \frac{1}{2}}{\frac{\sigma}{2}}\right)^2}\mathrm{d \left(\alpha_0(K)^{-1}(\nu_0^1)(t)\right)}.
$$
Since $\nu$ is a probability measure on $(0,\infty)$ and $\alpha_0(K)^{-1}(\nu_0^1)$ has uniform distribution on $[\frac{1}{2},1)$ so we have that
$$
p(\mathcal{K}_0)=\frac{1}{2\sqrt{2\pi}\frac{\sigma}{2}}\frac{1}{2}\left(\int\limits_{-\infty}^{1}e^{-\frac{1}{2}\left(\frac{t- \frac{1}{2}}{\frac{\sigma}{2}}\right)^2}\mathrm{dt}-\int\limits_{-\infty}^{\frac{1}{2}}e^{-\frac{1}{2}\left(\frac{t- \frac{1}{2}}{\frac{\sigma}{2}}\right)^2}\mathrm{dt}\right)=
$$
$$
=\frac{\Phi\left(\frac{1}{\sigma}\right)-\Phi(0)}{4},
$$
where the function
$$
\Phi(x) = \frac{1}{\sqrt{2\pi}} \int_\infty^x e^{\left(-\frac{u^2}{2}\right)} \mathrm{ du}
$$
is the distribution function of the standard normal distribution.

Analogously, for the set $\mathcal{K}_0(c):=\left\{K'\in \mathcal{K}_{0} \mid \alpha_0(K')=\alpha_0(K)\leq c\right\}$ we have
$$
p(\mathcal{K}_0(c))=\frac{\nu((0,\infty))}{\sqrt{2\pi\sigma^2}}\int\limits_{\frac{1}{2}}^{c}e^{-\frac{1}{2}\left(\frac{t- \frac{1}{2}}{\frac{\sigma}{2}}\right)^2}\mathrm{d \left(\alpha_0(K)^{-1}(\nu_0^1)(t)\right)}
=\frac{\Phi\left(\frac{c-\frac{1}{2}}{\frac{\sigma}{2}}\right)-\Phi(0)}{4}
$$
thus the measure
$$
P(\mathcal{A}):=\frac{4}{\Phi\left(\frac{1}{\sigma}\right)-\Phi(0)}p(\mathcal{A})
$$
is such a probability measure on $\mathcal{K}_0$ which pushforward by the function $\alpha_0(K)$ has normal distribution.
\qed
\begin{corollary}
From Theorem 2 follows the existence of a measure with similar properties on the space $\mathcal{K}$ of convex bodies. Let denote by  $m(K):=\frac{1}{2}(K+(-K))$ where the addition means the Minkowski sum of convex bodies. The mapping
$$
m:\mathcal{K}\longrightarrow \mathcal{K}_0
$$
is a continuous function on $\mathcal{K}$ and thus it defines a pullback measure $\mu$ on $\mathcal{K}$ by the rule
$$
\mu(H)=P(m(H)) \mbox{ where } H=m^{-1}(H') \mbox{ for a Borel set } H'\in \mathcal{K}_0.
$$
Observe that $m$ has the following properties:
\begin{enumerate}
\item surjective
\item for any set $S\subset \mathcal{K}$ and a vector $t\in \mathbb{R}^n$ we have $m(S+t)=m(S)$
\item for any $K\in \mathcal{K}$ holds that $d(K)=d(m(K))$, $w(K)=w(m(K)$ implying that $\alpha_0(K)=\alpha_0(m(K))$.
\end{enumerate}
This implies that the function $\alpha_0$ is well-defined on $\mathcal{K}$ and for any Borel set $B\in \left[\frac{1}{2},1\right)$
$$
\mu\left(\alpha_0^{-1}(B)\right)=P(m\left(\alpha_0^{-1}(B)\right))=P(\alpha_0^{-1}|_{\mathcal{K}_0}(B)),
$$
showing that the pushforward of the measure $\mu$ has truncated normal distribution on the interval $\left[\frac{1}{2},1\right)$.
\end{corollary}

\section{Geometric measure with normal pushforward.}

In this section we reformulate the preceding construction such a way it will be useful to stochastic-geometric examination. The basic questions on such a measure are the followings: " do the convex polytopes have measure zero,
do the smooth bodies have positive measure, or does a neighbourhood always have positive measure?" The improved construction gives positive answer to these questions. Since there is no polytope among the bodies $K_{m(n,k)}$ it seems to be immediately that the measure of the set of $O$-symmetric polytopes $\mathcal{P}_0$ are equal to zero. In fact, we introduce the sets $\mathcal{P}_0^1$ and $\widetilde{\mathcal{P}_0^1}$ as we did in the case of the $O$-symmetric bodies $\mathcal{K}_0$. By definition
$$
\widetilde{\mu_0^1}\left(\widetilde{\mathcal{K}_0^1}\setminus \widetilde{\mathcal{P}_0^1}\right)=1,
$$
showing that
$$
\widetilde{\mu_0^1}\left(\widetilde{\mathcal{P}_0^1}\right)=0.
$$
Thus
$$
\widetilde{\nu_0^1}\left(\widetilde{\mathcal{P}_0^1}\right)=\int\limits_{\widetilde{\mathcal{P}_0^1}}\mathrm{d}\widetilde{\nu_0^1}= \int\limits_{\widetilde{\mathcal{P}_0^1}}\frac{4}{\left(w+4\right)^2}\mathrm{d}\widetilde{\mu_0^1}=0,
$$
and so
$$
\nu_0^1\left(\mathcal{P}_0^1\right)=\Phi_{2a}^1\left(\widetilde{\mathcal{P}_0^1}\otimes \nu_n\right)\left(\widetilde{\mathcal{P}_0^1},SO(n)\right)=0.
$$
Finally, we have
$$
\nu_0\left(\mathcal{P}_0\right)=\Phi_4\left(\nu_0^1\otimes \nu\right)\left(\mathcal{P}_0^1,[0,\infty)\right)=0,
$$
hence
$$
p\left(\mathcal{P}_0\right)=P\left(\mathcal{P}_0\right)=0
$$
as we stated. It is clear that this calculation for the set of smooth $O$-symmetric bodies gives measure $1$, if we change the bodies $K_{m(n,k)}$ to smooth bodies. On the other hand the measure of a neighborhood will be positive, if and only if the system of bodies for which the measure concentrated will be dense in $\widetilde{\mathcal{K}_0^1}$ with respect to the Hausdorff metric. Thus the improvement of the construction change the system $K_{m(n,k)}$ to such a system of bodies of $\widetilde{\mathcal{K}_0^1}$.

{\bf In the first step} we change the body $K_{m(n,k)}$ to a smooth body $K_{m(n,k)}^l$ defined by the convex hull of the ball around the origin with radius $m(n,k)$ and the two balls of radius $\varepsilon_l=\frac{1}{2^l}m(n,k)$ with centers $\pm (2-\varepsilon_l)e_n$.

{\bf In the second step} we consider a countable system of bodies substituting the body $K_{m(n,k)}^l$ with the property that a given polytope $Q$ and a given number $\varepsilon$ we can choose an element $R$ from the system with $\delta^h\left(Q, R\right)<\varepsilon$. Consider a dense countable and centrally symmetric point system $\{ P_1,-P_1,P_2,-P_2\cdots \}$ in the closed ball of radius $2$ with the additional property that there is no two distances between the pairs of points which are equals to each other. (Such a point system is exist.) We assume that the first point $P_1$ is the endpoint of $2e_n$ and denote by $S_i$ a similarity of $E^n$ which sends $P_1$ into $P_i$ and the ball of radius $2$ at the origin into the ball of radius $OP_i$ centered at the origin $O$, too. Consider the countable set of bodies
$$
S\left(K_{m(n,k)}^l\right):=\left\{ S_i\left(K_{m(n,k)}^l\right) \mbox{ , } i=1,2,\ldots \right\},
$$
and define the elements of the set $\mathcal{H}_{m(n,k)}^l$ as follows:
\begin{itemize}
\item The first element is $K_{m(n,k)}^l:=S_1\left(K_{m(n,k)}^l\right)$
\item In the second step consider such pairs from the list of bodies one of which has diameter $4$ and add their convex hulls.
\item In the third step add the convex hull of the triplet from which one has dimension $4$.
\item ... and so on.
\end{itemize}
Hence we have a countable system of centrally symmetric bodies with diameter $4$, which partitioned (by the number of bodies from which an element generated as them convex hull) into countable subsets. So we have
$$
\mathcal{H}_{m(n,k)}^l=K_{m(n,k)}^l\cup \left\{ \mathrm{ conv }\left\{ S_i\left(K_{m(n,k)}^l\right), S_j\left(K_{m(n,k)}^l\right)\right\} \mbox{ for } i,j\right\} \cup
$$
$$
\cup \left\{ \mathrm{ conv }\left\{ S_i\left(K_{m(n,k)}^l\right), S_j\left(K_{m(n,k)}^l\right), S_k\left(K_{m(n,k)}^l\right)\right\}\mbox{ for } i,j,k \right\} \cup \cdots
$$
Note that the bodies of $\mathcal{H}_{m(n,k)}^l$ pairwise non-congruent. We distribute that part of the measure $\widetilde{\mu_0^1}$ which originally concentrated on $K_{m(n,k)}^l$ among the elements of $\mathcal{H}_{m(n,k)}^l$. If $L_i^r(l)$ is the $i^{th}$ element of the subset containing the convex hull of exactly $r$ copies of $S\left(K_{m(n,k)}^l\right)$ we add for it the weight $\frac{\alpha _i^r}{2^r}$, where for a fixed $r$ the sequences $(\alpha_i^r)$ of positive numbers holds the properties $\sum\limits_{i=1}^{\infty}\alpha_i^r=1$. Finally, choose a sequence of positive numbers $\beta_l$ with property $\sum\limits_{l=1}^{\infty}\beta_l=1$ and define the measure $\widetilde{\mu_0^1}$ as follows
$$
\widetilde{\mu_0^1}:= \lim\limits_{n\rightarrow \infty}\sum \limits_{k=1}^{2^{n+1}}\sum\limits_{l=1}^{\infty}\sum\limits_{r=1}^{\infty}\sum\limits_{i=1}^{\infty}\frac{\beta_l \alpha_i^r}{2^{n+1+r}}\delta_{L_i^r(l)}.
$$
With this measure all of our preceding statements are valid. It is easy to see that the required property on the approximation of polytopes is also hold because of the fact that for large $l$, $m(n,k)$ with a small $k$ the bodies $S\left(K_{m(n,k)}^l\right)$ essentially are $O$-symmetric segments, and thus for each polytope we can find a body from $\mathcal{H}_{m(n,k)}^l$ close to them. On the other hand the observations at the beginning of this section can be used thus the  $P$-measure of the set of polytopes is zero and the $P$-measure  of the set of smooth bodies is $1$, respectively.
We thus proved the following theorem:

\begin{theorem}
On the space of norms there is a probability measure $P$ with the following properties:
\begin{itemize}
\item The neighborhoods has positive measure.
\item The set of polytopes has zero measure.
\item The set of smooth bodies has measure 1.
\item The pushforward $\alpha_0(K)^{-1}(P)$ of $P$ has truncated normal distribution on the interval $[\frac{1}{2},1)$.
\end{itemize}
\end{theorem}

\begin{center}
\'A.G.Horv\'ath\\ Department of Geometry, Mathematical Institute\\
Budapest University of Technology and Economics,\\
H-1521 Budapest,\\
Hungary \\
e-mail: ghorvath@math.bme.hu
\end{center}

\end{document}